\documentclass[10pt]{article}
\usepackage{amsfonts}
\usepackage{graphicx}
\usepackage{latexsym}
\usepackage{amssymb}
\usepackage{amsmath}
\usepackage{layout}

\newtheorem{prop}{Proposition}
\newtheorem{lemma}{Lemma}

\newtheorem{corollary}{Corollary}
\newtheorem{theorem}{Theorem}
\newtheorem{remark}{Remark}

\def\real{{\mathord{{\rm I\kern-2.8pt R}}}}
\def\inte{{\mathord{{\rm I\kern-2.8pt N}}}}

\def\sZZ{{\rm Z\kern-2.8ptem{}Z}}

\def\z{{\mathchoice
  {\sZZ}
  {\sZZ}
  {\rm Z\kern-0.30em{}Z}
  {\rm Z\kern-0.25em{}Z} }}
\def\sQQ{{\kern 0.27em \vrule height1.45ex width0.03em depth0em
          \kern-0.30em \rm Q}}
\def\qu{{\mathchoice
    {\sQQ}
    {\sQQ}
  {\kern 0.225em \vrule height1.05ex width0.025em depth0em \kern-0.25em \rm Q}
  {\kern 0.180em \vrule height0.78ex width0.020em depth0em \kern-0.20em \rm Q}
        }}
\def\sCC{{\kern 0.27em \vrule height1.45ex width0.03em depth0em
          \kern-0.30em \rm C}}
\def\complex{{\mathchoice
    {\sCC}
    {\sCC}
  {\kern 0.225em \vrule height1.05ex width0.025em depth0em \kern-0.25em \rm C}
  {\kern 0.180em \vrule height0.78ex width0.020em depth0em \kern-0.20em \rm C}
        }}

\newcommand{\ba}{\begin{array}}
\newcommand{\ea}{\end{array}}
\newcommand{\be}{\begin{equation}}
\newcommand{\ee}{\end{equation}}
\newcommand{\bea}{\begin{eqnarray}}
\newcommand{\eea}{\end{eqnarray}}
\newcommand{\beaa}{\begin{eqnarray*}}
\newcommand{\eeaa}{\end{eqnarray*}}

\def\z{\zeta}

\font\tenmath=msbm10 \font\sevenmath=msbm7 \font\fivemath=msbm5
\newfam\mathfam \textfont\mathfam=\tenmath
\scriptfont\mathfam=\sevenmath \scriptscriptfont\mathfam=\fivemath

\def \={{\buildrel {\rm (law)} \over =}}

\newcommand{\basa}{\begin{assumption}}
\newcommand{\easa}{\end{assumption}}

\newcommand{\bas}{\begin{assum}}
\newcommand{\eas}{\end{assum}}

\newcommand{\ignore}[1]{}
\begin{document}

\title{Anticipating integrals and martingales on the Poisson space}
\author{Giovanni Peccati $^{1}\quad$Ciprian A. Tudor $^{2}\vspace*{0.1in}$ \\
$^{1}$ Laboratoire de Statistique Th\'eorique et Appliqu\'ee, Universit\'{e}
de Paris 6,\\
4, place Jussieu, F-75252 Paris Cedex 05, France.\\
giovanni.peccati@gmail.com\\
$^{2}$SAMOS/MATISSE, Universit\'e de Panth\'eon-Sorbonne Paris 1,\\
90, rue de Tolbiac, F-75634 Paris Cedex 13, France.\\
tudor@univ-paris1.fr\vspace*{0.1in}}
\maketitle

\begin{abstract}
Let $\tilde{N}_{t}$ be a standard compensated Poisson process on
$[0,1]$. We prove a new characterization of anticipating integrals
of the Skorohod type with respect to $\tilde{N}$, and use it to
obtain several counterparts to well established properties of
semimartingale stochastic integrals. In particular we show that,
if the integrand is sufficiently regular, anticipating Skorohod
integral processes with respect to $\tilde{N}$ admit a pointwise
representation as usual It\^{o} integrals in an independently
enlarged filtration. We apply such a result to: (i) characterize
Skorohod integral processes in terms of products of backward and
forward Poisson martingales, (ii) develop a new It\^{o}-type
calculus for anticipating integrals on the Poisson space, and
(iii) write Burkholder-type inequalities for Skorohod integrals.


\end{abstract}

\section{Introduction}

Let $\tilde{N}_{t}$ be a standard compensated Poisson process on $\left[ 0,1%
\right] $. The aim of this paper is to prove a new characterization of
anticipating integrals (of the Skorohod type) with respect to $\tilde{N}$,
and to apply such a result to investigate the relations between anticipating
integrals and Poissonian martingales.

\bigskip

The anticipating Skorohod integral has been first introduced in \cite{Sko}
in the context of Gaussian processes. It is well known that the notion of
Skorohod integral can be naturally extended to the family of \textit{normal
martingales}, that is, martingales having a conditional quadratic variation
equal to $t$ (among which there are the Wiener process and the compensated
Poisson process). In this case, the Skorohod integral is an extension of the
classical (semimartingale) It\^{o} integral to a wider family of non-adapted
integrands, and therefore coincides with the latter on the class of (square
integrable) adapted processes. See, for instance, \cite{MPS}; we also refer
to \cite{N} for an exhaustive presentation of results, techniques and
applications of the anticipating stochastic calculus in the Gaussian
context. In this paper, we try to deal with some of the disadvantages of the
Skorohod calculus, in the specific case of the compensated Poisson process,
and from the standpoint of the classic semimartingale theory.

\bigskip

To better understand our motivations, consider a Skorohod integral process
(see \cite{MPS} for the precise setup) with the form%
\begin{equation}
X_{t}=\int_{0}^{t}u_{s}\delta Z_{s}\hskip0.7cmt\in \lbrack 0,1]  \label{unu}
\end{equation}%
where $Z$ is a normal martingale and $u\mathbf{1}_{[0,t]}$ belongs
to the domain of the Skorohod integral $\delta$, $\forall t$.
Then, in general, if the integrand $u$ is not adapted to the
natural filtration of $Z$, the process $X$ is not a
semimartingale, and the study of $X$ cannot be carried out by
means of the usual It\^{o} theory (as presented for instance in
\cite{Elliott}). Actually, the techniques employed to deal with
processes such as $X$ mostly stem from functional analysis, and
they do not allow, e.g., to obtain fine trajectorial properties.
However, when $Z=W$, where $W$ stands for a standard Wiener
process, the authors of the present paper have pointed out several
remarkable connections between processes such as $X$ in
(\ref{unu}) and Wiener martingales. In particular, the following
results (among others) have been obtained when $Z=W$:

\begin{description}
\item[i)] the class of Skorohod integral processes with sufficiently regular
integrands coincides with a special family of It\^{o} integrals, called
\textit{It\^{o}-Skorohod integrals} (see \cite{Tudor});

\item[ii)] a Skorohod integral process $X$ such as (\ref{unu}) can be
approximated, in a certain norm, by linear combinations of processes with
the form $M_{t}\times M_{t}^{\prime }$, where $M_{t}$ is (centered) Wiener
martingale and $M_{t}^{\prime }$ is a \textit{backward} Wiener martingale
(see \cite{PTT}).
\end{description}

\bigskip

These facts lead in a natural way to explore the anticipating integrals in
the context of a standard Poisson process. It is known that a Skorohod type
integration can be developed on the Poisson space by using the Fock space
structure, and that such integrals enjoy a number of useful properties, in
part analogous to the ones displayed by Skorohod operators on Wiener space.
We refer e.g. to \cite{NV}, \cite{Ocone}, \cite{CP}, \cite{LT}, \cite{P1},
\cite{P2} or \cite{MPS} for different aspects of this calculus. Here, we
shall provide Poissonian counterparts to several results given in \cite%
{Tudor} and \cite{PTT}, as facts (i) and (ii) above, and we shall
systematically point out the arguments that differ from those
given in the Wiener context. We remark that some of our results
are of particular
interest in the Poisson case. For example, our methods allow to obtain an It%
\^{o}-type formula for anticipating integrals, and -- as far as we
know -- this is the only anticipating change of variables formula
for the Poisson situation (we could not find, for instance, an
It\^{o}-type formula proved in the spirit of \cite{NP}). We note
that the fact that the increments of the Poisson process are
independent plays a crucial role in our construction; therefore the
extension of the results to a more general normal martingale seems
difficult.

\bigskip

The paper is organized as follows. The first part of Section 2
contains some preliminaries on the Malliavin calculus with respect
to the Poisson process, whereas the second part displays a
discussion about It\^{o} stochastic integrals, $\sigma $-fields
and (independently enlarged) filtrations on the Poisson space; we
shall note that most of the results given here are still valid on
the Wiener space. In Section 3, we show that every anticipating
(Poisson) integral coincides pointwise with a special type of
It\^{o} integral and (as in \cite{PTT}) we apply this relation to
approximate Skorohod integral processes by linear combinations of
processes that are a product of forward and backward martingales.
Finally, in Section 4 we develop a new stochastic calculus of the
It\^{o} type for anticipating Poisson integrals.

\section{Preliminaries: Malliavin calculus and filtrations}

\subsection{Malliavin calculus for the Poisson process}

Throughout the paper, we use notation and terminology from standard
semimartingale theory. The reader is referred e.g. to \cite{Elliott}, \cite%
{DMM} or \cite{Protter} for any unexplained notion.

\bigskip

Let $T=[0,1]$ and let $N=(N_{t})_{t\in T}$ be a standard Poisson process,
defined on the standard Poisson space $\left( \Omega ,\mathbb{F},\mathbf{P}%
\right) $ (see e.g. \cite{NV}). By $\tilde{N}$ we will denote the
compensated Poisson process $\tilde{N}_{t}=N_{t}-t$. For every Borel set $B$%
, we will note $N\left( B\right) $ and $\tilde{N}\left( B\right) $,
respectively, the random measures $\sum_{s\in B}\Delta N_{s}$ and $N\left(
B\right) -\lambda \left( B\right) $, where $\Delta N_{s}=N_{s}-N_{s-}$ and $%
\lambda $ stands for Lebesgue measure on $T$. It is well known (see e.g.
\cite{DMM} or \cite{MPS}) that the process $t\mapsto \tilde{N}_{t}$ is a
\emph{normal martingale}, that is, $\tilde{N}_{t}$ is a c\`{a}dl\`{a}g
martingale initialized at zero, such that its \emph{conditional quadratic
variation} process (or \emph{angle bracket} process) is given by $\langle
\tilde{N},\tilde{N}\rangle _{t}=\langle \tilde{N}\rangle _{t}=t$. The \emph{%
quadratic variation process} of $\tilde{N}$ (or \emph{right bracket process}%
) is of course $[\tilde{N},\tilde{N}]_{t}=[\tilde{N}]_{t}=N_{t}$ (note that
the results of our paper extend immediately to the case of a Poisson process
on $\mathbb{R}_{+}$, with a deterministic intensity $\mu >0$). It is also
known (see again \cite{MPS} and the references therein) that $\tilde{N}$
enjoys the \emph{chaotic representation property, }i.e. every random
variable $F\in L^{2}\left( \Omega ,\mathbb{F},\mathbf{P}\right) =L^{2}\left(
\mathbf{P}\right) $, measurable with respect to the $\sigma $-algebra
generated by $N$, can be written as an orthogonal sum of \emph{multiple
Poisson-It\^{o} stochastic integrals }
\begin{equation}
F=\mathbf{E}\left( F\right) +\sum_{n\geq 1}I_{n}(f_{n})  \label{chaos}
\end{equation}%
where the infinite series converges in $L^{2}\left( \mathbf{P}\right) $,
and, for $n\geq 1$, the kernel $f_{n}$ is an element of $L_{s}^{2}(T^{n})%
\subset L^{2}(T^{n})$, where $L_{s}^{2}(T^{n})$ and $L^{2}(T^{n})$ denote,
respectively, the space of symmetric and square integrable functions, and
the space of square integrable functions on $T^{n}$ (endowed with Lebesgue
measure).

\bigskip

Let us recall the basic construction of the multiple Poisson-It\^{o}
integral on the Poisson space. Fix $n\geq 2$ and denote by $S_{n}$ and $%
\widetilde{S}_{n}$, respectively, the vector space generated by simple
functions with the form%
\begin{equation}
f(t_{1},\ldots ,t_{n})=\mathbf{1}_{B_{1}}(t_{1})\ldots \mathbf{1}%
_{B_{n}}(t_{n})\text{,}  \label{sn}
\end{equation}%
where $B_{1},\ldots ,B_{n}$ are disjoint subsets of $\mathbb{R}$, and the
vector space generated by the symmetrization of the element of $S_{n}$. If $%
f $ is as in (\ref{sn}) and $\tilde{f}\in \widetilde{S}_{n}$ is its
symmetrization, we define $I_{n}(\tilde{f})$ as
\begin{equation}
I_{n}(\tilde{f})=\tilde{N}(B_{1})\ldots \tilde{N}(B_{n})  \label{in}
\end{equation}%
Since, for every $n\geq 2$, $\widetilde{S}_{n}$ is dense in $%
L_{s}^{2}(T^{n}) $, the integral $I_{n}$ can be extended to $%
L_{s}^{2}(T^{n}) $ by continuity, due to the isometry formula, true for
every $m,n\geq 2$, $\tilde{f}\in \widetilde{S}_{n}$ and $\tilde{g}\in
\widetilde{S}_{m},$%
\begin{equation}
\mathbf{E}\left( I_{n}(\tilde{f})I_{m}(\tilde{g})\right) =n!\langle \tilde{f}%
,\tilde{g}\rangle _{L^{2}(T^{n})}\mathbf{1}_{(n=m)}.  \label{iso}
\end{equation}%
We also use the following conventional notation: $L^{2}\left( T\right) =$ $%
L^{2}\left( T^{1}\right) =L_{s}^{2}\left( T^{1}\right) $; $I_{1}\left(
f\right) =\int_{0}^{1}f\left( s\right) d\tilde{N}_{s}$, $f\in L^{2}\left(
T\right) $; $\tilde{f}$ is the symmetrization of $f$, $\forall f\in
L^{2}\left( T^{n}\right) $, $n\geq 2$; $I_{n}\left( f\right) =I_{n}( \tilde{f%
}) $, $f\in L^{2}\left( T^{n}\right) $, $n\geq 2$; $L^{2}\left( T^{0}\right)
=L_{s}^{2}\left( T^{0}\right) =$ $S_{0}=$ $\widetilde{S}_{0}=\mathbb{R}$; $%
I_{0}\left( c\right) =c$, $c\in \mathbb{R}$.

\bigskip

\begin{remark}
As proved e.g. by Ogura in \cite{Ogura}, one can define multiple stochastic
integrals on the Poisson space by using the Charlier-Poisson polynomials.
More precisely, for $n\geq 0$, the $n$th Charlier-Poisson polynomial $%
C_{n}(t,x)$, $\left( t,x\right) \in \left[ 0,1\right] \times \mathbb{R}$, is
defined through the generating function (see for instance \cite{AS})
\begin{equation*}
\Phi (z,t,x)=\sum_{n=0}^{\infty }z^{n}C_{n}(t,x)=(1+z)^{x+t}\exp (-zt).
\end{equation*}%
It is well known (see e.g. \cite[Lemma 2]{Kab}) that the Charlier
polynomials are connected to the above defined Poisson-It\^{o} multiple
integrals by the following relation: for every Borel subset $B\subseteq T$
\begin{equation}
C_{n}(\lambda \left( B\right) ,\tilde{N}\left( B\right) )=\frac{1}{n!}%
I_{n}\left( 1_{B}^{\otimes n}(\cdot )\right) ,  \label{CharlierP}
\end{equation}%
where $n\geq 1$, and $\lambda $ stands for Lebesgue measure.
\end{remark}

\bigskip

Now define, for $n,m\geq 1$, $f\in L_{s}^{2}(T^{m})$, $g\in L_{s}^{2}(T^{n})$%
, $r=0,...,m\wedge n$ and $l=1,...,r$, the (contraction) kernel on $%
T^{m+n-r-l}$%
\begin{eqnarray*}
&&f\star _{r}^{l}(\gamma _{1},\ldots ,\gamma _{r-l},t_{1},\ldots
,t_{m-r},s_{1},\ldots ,s_{n-r}) \\
&=&\int_{T^{l}}f(u_{1},\ldots ,u_{l},\gamma _{1},\ldots ,\gamma
_{r-l},t_{1},\ldots ,t_{m-r})g(u_{1},\ldots ,u_{l},\gamma _{1},\ldots
,\gamma _{r-l},s_{1},\ldots ,s_{n-r})du_{1}...du_{l}\text{,}
\end{eqnarray*}%
and, for $l=0$,
\begin{equation*}
f\star _{r}^{0}(\gamma _{1},\ldots ,\gamma _{r},t_{1},\ldots
,t_{m-r},s_{1},\ldots ,s_{n-r})=f(\gamma _{1},\ldots ,\gamma
_{r},t_{1},\ldots ,t_{m-r})g(\gamma _{1},\ldots ,\gamma _{r},s_{1},\ldots
,s_{n-r}),
\end{equation*}%
so that
\begin{equation*}
f\star _{0}^{0}(t_{1},\ldots ,t_{m},s_{1},\ldots ,s_{n})=f(t_{1},\ldots
,t_{m})g(s_{1},\ldots ,s_{n}).
\end{equation*}
We will need the following product formula for two Poisson-It\^{o} multiple
integrals (see \cite{Kab}, \cite{P3}, or \cite{CT}): let $f\in
L_{s}^{2}(T^{m})$ and $g\in L_{s}^{2}(T^{n})$, $n,m\geq 1$, and suppose
moreover that $f\star _{r}^{l}g\in L^{2}(T^{m+n-r-l})$ for every $%
r=0,...,m\wedge n$ and $l=1,...,r$, then
\begin{equation}
I_{m}(f)I_{n}(g)=\sum_{r=0}^{m\wedge n}r!\dbinom{m}{r}\dbinom{n}{r}%
\sum_{l=0}^{r}I_{m+n-r-l}(f\star _{r}^{l}g).  \label{product}
\end{equation}

\bigskip

It is possible to develop a Malliavin-type calculus with respect to the
Poisson process based on the (symmetric) Fock space isomorphism induced by
formulae (\ref{chaos}) and (\ref{iso}). We refer to \cite{NV} or \cite{MPS}
for the basic elements of this calculus. For a random variable $F$ as in (%
\ref{chaos}) we introduce the \emph{annihilation} (or \emph{derivative})
operator as
\begin{equation}
D_{t}F=\sum_{n\geq 1}nI_{n-1}(f_{n}(\cdot ,t)),\hskip0.5cmt\in T ,
\label{deriv}
\end{equation}%
and its domain, usually denoted by $\mathbb{D}^{1,2}$, is the set
\begin{equation*}
\mathbb{D}^{1,2}=\{F=\sum_{n\geq 0}I_{n}(f_{n}):\sum_{n}nn!\Vert f_{n}\Vert
_{n}^{2}<+\infty \}
\end{equation*}%
where $\Vert \cdot \Vert _{n}$ is the norm in $L^{2}(T^{n})$. The operator $%
D $ is not a derivation (see \cite[p. 91]{MPS}) and it satisfies (see \cite[%
Lemma 6.1\ and Theorem 6.2]{NV})
\begin{equation}
D(FG)=FDG+GDF+DFDG\mbox{ if }F,G,FG\in \mathbb{D}^{1,2}.  \label{chain rule}
\end{equation}

\bigskip

The \emph{Skorohod integral}, or the \emph{creation} operator, is defined by
\begin{equation*}
\delta (u)=\sum_{n\geq 0}I_{n+1}(\tilde{f}_{n})
\end{equation*}%
whenever $u_{t}=\sum_{n\geq 0}I_{n}(f_{n}(\cdot ,t))$, where $u\in
L^{2}(T\times \Omega )$, belongs to the domain of $\delta $, noted $%
Dom(\delta )$, that is, $u$ verifies
\begin{equation*}
\sum_{n\geq 0}nn!\Vert f_{n}\Vert _{n+1}^{2}<+\infty .
\end{equation*}

We introduce the subset $\mathbb{L}^{1,2}$ of $Dom\left( \delta \right) $
defined as
\begin{equation}
\mathbb{L}^{1,2}=\{u_{t}=\sum_{n}I_{n}(f_{n}(\cdot ,t)):\sum_{n}(n+1)!\Vert
f_{n}\Vert _{n+1}^{2}<+\infty \}.  \label{l12}
\end{equation}%
Note that $\mathbb{L}^{1,2}$ equals $L^{2}(T;\mathbb{D}^{1,2})$, when the
former is endowed with the seminorm
\begin{equation*}
\Vert u\Vert _{1,2}^{2}=\mathbf{E}\int_{0}^{1}u_{s}^{2}ds+\mathbf{E}%
\int_{0}^{1}\int_{0}^{1}(D_{r}u_{s})^{2}drds,
\end{equation*}%
and moreover, for every $u\in \mathbb{L}^{1,2}$, one can verify the
inequality
\begin{equation}
\mathbf{E}\delta (u)^{2}\leq \Vert u\Vert _{1,2}^{2}  \label{meyer}
\end{equation}%
(see for instance formula (\ref{cov}) below). For any integer $k\geq 2$ the
space $\mathbb{D}^{k,2}$ denotes the set of $k$ times \textit{weakly
differentiable} random variables, endowed with the seminorm
\begin{equation*}
\Vert F\Vert _{k,2}^{2}=\mathbf{E}\Vert F\Vert ^{2}+\sum_{l=1}^{k}\Vert
D^{l}F\Vert _{L^{2}(T^{l}\times \Omega )}^{2}, \qquad F\in\mathbb{D}^{k,2},
\end{equation*}%
where $D^{1}=D$, and, for $l\geq 2$, the $l$th Malliavin derivative $D^{l}$
is first defined by iteration on simple functionals, and then extended by a
standard closure procedure. By $\mathbb{L}^{k,2}$ we denote the set $L^{2}(T;%
\mathbb{D}^{k,2})$; note that, if $u\in \mathbb{L}^{1,2}$, then $u\mathbf{1}%
_{\left[ 0,t\right] }\in Dom\left( \delta \right) $ for every $t\in T$ (see
\cite[Lemma 4.6]{MPS}). Occasionally, we will also use the notation%
\begin{equation*}
\int_{0}^{1}u_{s}\delta \tilde{N}_{s}=\delta \left( u\right) \text{, \ \ }%
u\in Dom\left( \delta \right) \text{.}
\end{equation*}%
We recall that, according e.g. to \cite[Proposition 4.4]{MPS}, if $u\in
L^{2}\left( T\times \Omega \right) $ is a predictable process with respect
to the filtration generated by $\tilde{N}$, then $u$ is Skorohod integrable
and $\delta \left( u\right) $ coincides with the usual It\^{o} integral with
respect to the c\`{a}dl\`{a}g martingale $\tilde{N}$.

\bigskip

The following duality relationship between $D$ and $\delta $, which is
classic in the Gaussian case, still holds on the Poisson space (see \cite[%
Proposition 4.2]{NV})%
\begin{equation}
\mathbf{E}\left[ \int_{0}^{1}\left( D_{s}Fu_{s}\right) ds\right] =\mathbf{E}%
\mathbb{(}\delta (u)F)\mbox{ if }F\in \mathbb{D}^{1,2}\mbox{ and } u\in
Dom(\delta ).  \label{dua}
\end{equation}%
and the covariance of two Skorohod integrals is given by
\begin{equation}
\mathbf{E}(\delta (u)\delta (v))=\mathbf{E}\int_{0}^{1}u_{s}v_{s}ds+\mathbf{E%
}\int_{0}^{1}\int_{0}^{1}D_{r}u_{s}D_{s}v_{r}drds  \label{cov}
\end{equation}%
whenever $u,v\in \mathbb{L}^{1,2}$ (see \cite[Thorem 4.1]{NV}). Eventually,
we will need the commutativity relationship between $D$ and $\delta $
\begin{equation}
D_{t}\delta (u)=u_{t}+\delta (D_{t}u)\mbox{ if }u\in \mathbb{L}^{1,2}%
\mbox{
and }D_{t}u\in Dom(\delta )  \label{comm}
\end{equation}%
(see \cite[Thorem 4.1]{NV}), as well as the integration by parts formula,
true whenever $F\in \mathbb{D}^{1,2}\ $and $DF\times u\in Dom\left( \delta
\right) ,$%
\begin{equation}
\delta (Fu)=F\delta (u)-\langle DF,u\rangle _{L^{2}(T)}-\delta (DFu)
\label{ip}
\end{equation}%
(see \cite[Theorem 7.1]{NV}).

\subsection{$\protect\sigma $-fields, filtrations and It\^{o} stochastic
integrals}

For any Borel set $A\subseteq T$, we use the notation
\begin{equation*}
\mathbb{F}_{A}=\sigma \left\{ \tilde{N}\left( h\mathbf{1}_{A}\right) :h\in
L^{2}\left( T\right) \right\}
\end{equation*}%
and also, for $t\in T$, $\mathbb{F}_{t}=\mathbb{F}_{\left[ 0,t\right] }$ and
$\mathbb{F}_{t^{c}}=\mathbb{F}_{\left[ 0,t\right] ^{c}}$. Note that, in the
following, we will tacitly complete each $\sigma $-field $\mathbb{F}_{A}$
with the $\mathbf{P}$-null sets of $\mathbb{F}$, so that, for instance, $%
\left\{ \mathbb{F}_{t}:t\in T\right\} $ is the completion of the natural
filtration of the process $\tilde{N}_{t}$. We also set $\mathbb{F=F}_{1}$.
The independence of the Poisson increments implies the following relations
(see again \cite{NV}): for every $n\geq 1$, every $f\in L_{s}^{2}\left(
T^{n}\right) $ and every Borel subset $A$ of $T$,
\begin{equation}
\mathbf{E}\left( I_{n}(f)\mid \mathbb{F}_{A}\right) =I_{n}\left( f\mathbf{1}%
_{A}^{\otimes n}\right) \text{, a.s.-- }\mathbf{P}\text{,}  \label{r1}
\end{equation}%
and, a.s.-- $\mathbf{P}$,
\begin{equation}
D_{t}\mathbf{E}(F\mid \mathbb{F}_{A})=\mathbf{E}(D_{t}F\mid \mathbb{F}_{A})%
\mathbf{1}_{A}(t), \hskip0.5cmt\in T.  \label{r2}
\end{equation}%
An immediate consequence of (\ref{r2}) is that, if $F\in \mathbb{D}^{1,2}$
is a $\mathbb{F}_{A}$-measurable random variable, then $DF=0$ on $%
A^{c}\times \Omega .$ Moreover, if we denote by $X$ the Skorohod integral
process
\begin{equation*}
X_{t}=\int_{0}^{t}u_{s}\delta \tilde{N}_{s}=\delta (u\mathbf{1}_{[0,t]}),%
\hskip0.5cmt\in T,u\in \mathbb{L}^{1,2},
\end{equation*}%
then the process $X$ satisfies (see Lemma 3.2.1 in \cite{N})
\begin{equation}
\mathbf{E}\left( X_{t}-X_{s}\mid \mathbb{F}_{[s,t]^{c}}\right) =0%
\mbox{ for
every }s<t.  \label{r3}
\end{equation}
In the anticipating calculus, relation (\ref{r3}) plays roughly the same
role as does, in the usual It\^{o} calculus, the martingale characterization
of adapted stochastic integrals.

\bigskip

Now fix $t\in \left( 0,1\right] $. In the sequel, we will use the properties
of the following (enlarged) filtration
\begin{equation}
\mathbb{F}_{(\cdot ,t]^{c}}=\left\{ \mathbb{F}_{(s,t]^{c}}:s\in \left[ 0,t%
\right] \right\} =\left\{ \mathbb{F}_{s}\vee \mathbb{F}_{t^{c}}:s\in \left[
0,t\right] \right\} .  \label{filtraziun}
\end{equation}%
Note that, since $\tilde{N}$ has c\`{a}dl\`{a}g paths, $\mathbb{F}%
_{(s,t]^{c}}=\mathbb{F}_{(s,t)^{c}}$ for every $s\in \left[ 0,t\right] $,
and also
\begin{equation}
\mathbb{F}_{[s,t]^{c}}=\bigvee\limits_{\varepsilon >0}\mathbb{F}%
_{(s-\varepsilon ,t]^{c}}.  \label{Vfields}
\end{equation}%
It is also easily checked that the filtration $\mathbb{F}_{(\cdot ,t]^{c}}$
satisfies the usual conditions. We can therefore define, for every $t\in
\left( 0,1\right] $, $^{\left( p,t\right) }\left( \cdot \right) $ to be the
predictable projection operator with respect to $\mathbb{F}_{(\cdot ,t]^{c}}$%
, as defined e.g. in \cite[Theorem 6.39]{Elliott}. Note that, due to the
independence of its increments, the process $\tilde{N}_{r}$, $r\in \left[ 0,t%
\right] $, is again a normal martingale with respect to the filtration $%
\mathbb{F}_{(r,t]^{c}}$. It follows that the It\^{o} (semimartingale)
stochastic integral of a square integrable, $\mathbb{F}_{(\cdot ,t]^{c}}$ --
predictable process is always well defined. For a process $u\in L^{2}\left(
T\times \Omega \right) $, whose restriction to $\left[ 0,t\right] $ is also $%
\mathbb{F}_{(\cdot ,t]^{c}}$ -- predictable, we will note $\int_{0}^{t}u_{s}d%
\tilde{N}_{s}$ the It\^{o} stochastic integral of $u$ with respect to $%
\tilde{N}$, regarded as a c\`{a}dl\`{a}g, square integrable $\mathbb{F}%
_{(\cdot ,t]^{c}}$ -- martingale on $\left[ 0,t\right] $. Note that we write
$\int_{0}^{t}$ instead of $\int_{0+}^{t}$, because $\tilde{N}_{0}=0$. The
following result extends \cite[Proposition 4.4]{MPS} to the case of the
enlarged filtration $\mathbb{F}_{(\cdot ,t]^{c}}$. It also contains a
Clark-Ocone type formula (see e.g. \cite{Ocone} for the Brownian case) which
will be further generalized in the next section.

\begin{prop}
Let the above notation prevail and fix $t\in \left( 0,1\right] $. Then,

\begin{description}
\item[ ] (i) if the restriction to $\left[ 0,t\right] $ of a process $u\in
L^{2}\left( T\times \Omega \right) $ is predictable with respect to the
filtration $\mathbb{F}_{(\cdot ,t]^{c}}$, then
\begin{equation*}
u\mathbf{1}_{\left[ 0,t\right] }\in Dom\left( \delta \right) \text{ \ \ and
\ \ }\delta \left( u\mathbf{1}_{\left[ 0,t\right] }\right)
=\int_{0}^{t}u_{s}d\tilde{N}_{s},
\end{equation*}%
where the right-hand side is a stochastic integral in the semimartingale
sense;

\item[ ] (ii) for every $\mathbb{F}_{t}$-measurable functional $G\in \mathbb{%
D}^{1,2}$, and for every $0\leq s\leq t$,
\begin{equation}
G=\mathbf{E}\left( G\mid \mathbb{F}_{(s,t]^{c}}\right) +\delta \left(
^{\left( p,t\right) }\left( D_{\cdot }G\right) \right) =\mathbf{E}\left(
G\mid \mathbb{F}_{(s,t]^{c}}\right) +\int_{s}^{t}\text{ }^{\left( p,t\right)
}\left( D_{r}G\right) d\tilde{N}_{r}.  \label{MPS-CO}
\end{equation}
\end{description}
\end{prop}

\textbf{Proof. }(i) Fix $t\in \left( 0,1\right] $. We start by considering a
process with the form
\begin{equation}
u_{s}=g\left( s\right) I_{n}\left( h^{\otimes n}\mathbf{1}_{\left[
0,s\right) \cup \left( t,1\right] }^{\otimes n}\right) ,\text{ \ \ }s\in T%
\text{,}  \label{prodproc}
\end{equation}%
where $n\geq 1$, $h^{\otimes n}\left( t_{1},...,t_{n}\right) =h\left(
t_{1}\right) \cdot \cdot \cdot h\left( t_{n}\right) $, and $h,g\in
L^{2}\left( T\right) $. Plainly, $u\mathbf{1}_{\left[ 0,t\right] }\in
Dom\left( \delta \right) $. Now define, for $k=0,...,n$ and $s<t$%
\begin{equation*}
B_{k}^{n}\left( s\right) =\left\{ \left( t_{1},...,t_{n}\right) \in T^{n}:k%
\text{ of the }t_{i}\text{'s are }<s\text{ and }n-k\text{ are }>t\right\} .
\end{equation*}%
Of course,%
\begin{equation*}
u_{s}=\sum_{k=0}^{n}u_{s}^{k},\text{ \ \ }s\in T\text{,}
\end{equation*}%
where, for $k=0,...,n$,%
\begin{equation*}
u_{s}^{k}=g\left( s\right) I_{n}\left( h^{\otimes n}\mathbf{1}%
_{B_{k}^{n}\left( s\right) }\right) =\dbinom{n}{k}g\left( s\right)
I_{k}\left( h^{\otimes k}\mathbf{1}_{\left[ 0,s\right) }^{\otimes k}\right)
I_{n-k}\left( h^{\otimes n-k}\mathbf{1}_{\left( t,1\right] }^{\otimes
n-k}\right) ,
\end{equation*}%
with $I_{0}=1$, the last equality being justified by an application of (\ref%
{product}). Now fix $k$, and observe that the process $u_{s}^{k}$, $s\in %
\left[ 0,t\right] $, is predictable with respect to $\mathbb{F}_{(\cdot
,t]^{c}}$ so that, on $\left[ 0,t\right] $, the It\^{o} integral of $u^{k}$
with respect to $\tilde{N}$ is classically given by
\begin{equation*}
\int_{0}^{t}u_{s}^{k}d\tilde{N}_{s}=\dbinom{n}{k}I_{n-k}\left( h^{\otimes
n-k}\mathbf{1}_{\left( t,1\right] }^{\otimes n-k}\right) \int_{0}^{t}g\left(
s\right) I_{k}\left( h^{\otimes k}\mathbf{1}_{\left[ 0,s\right) }^{\otimes
k}\right) d\tilde{N}_{s}.
\end{equation*}%
On the other hand, the symmetrization in the $n+1$ variables $\left(
t_{1},t_{2},...,t_{n+1}\right) $ of the function $g\mathbf{1}_{\left[ 0,t%
\right] }\left( t_{1}\right) $ $\times $ $h^{\otimes n}\mathbf{1}%
_{B_{k}^{n}\left( t_{1}\right) }\left( t_{2},...,t_{n+1}\right) $ is
\begin{equation*}
f_{n+1}^{k}\left( t_{1},...,t_{n+1}\right) =\frac{1}{n+1}\sum_{i=1}^{n+1}g%
\mathbf{1}_{\left[ 0,t\right] }\left( t_{i}\right) h^{\otimes n}\mathbf{1}%
_{B_{k}^{n}\left( t_{i}\right) }\left( t_{j}:j\neq i\right) ,
\end{equation*}%
and the restriction of $f_{n+1}^{k}$ to $\Delta ^{n+1}=\left\{ \left(
t_{1},...,t_{n+1}\right) \in T^{n+1}:0<t_{1}<...<t_{n+1}<1\right\} $ is
therefore%
\begin{eqnarray*}
&&\frac{1}{n+1}g\mathbf{1}_{\left[ 0,t\right] }\left( t_{k+1}\right)
h^{\otimes n}\left( t_{1},...,t_{k},t_{k+2},...,t_{n+1}\right) \mathbf{1}%
_{\left\{ t_{k+2}>t\right\} } \\
&=&\frac{1}{n+1}g\mathbf{1}_{\left[ 0,t\right] }\left( t_{k+1}\right)
h^{\otimes k}\left( t_{1},...,t_{k}\right) h^{\otimes n-k}\left(
t_{k+2},...,t_{n+1}\right) \mathbf{1}_{\left\{ t_{k+2}>t\right\} }\text{, \
\ }\left( t_{1},...,t_{n+1}\right) \in \Delta ^{n+1}
\end{eqnarray*}%
and consequently
\begin{eqnarray*}
\delta \left( u^{k}\right) &=&n!\int_{t}^{1}\int_{t}^{t_{n+1}}\cdot \cdot
\cdot \int_{t}^{t_{k+3}}h^{\otimes n-k}\left( t_{k+2},...,t_{n+1}\right) d%
\tilde{N}_{t_{k+2}}d\tilde{N}_{t_{k+3}}...d\tilde{N}_{t_{n+1}}\times \\
&&\times \int_{0}^{t}g\left( t_{k+1}\right) \int_{0}^{t_{k+1}}\cdot \cdot
\cdot \int_{0}^{t_{2}}h^{\otimes k}\left( t_{1},...t_{k}\right) d\tilde{N}%
_{t_{1}}\cdot \cdot \cdot d\tilde{N}_{t_{k}} d\tilde{N}_{t_{k+1}} \\
&=&\dbinom{n}{k}I_{n-k}\left( h^{\otimes n-k}\mathbf{1}_{\left( t,1\right]
}^{\otimes n-k}\right) \int_{0}^{t}g\left( s\right) I_{k}\left( h^{\otimes k}%
\mathbf{1}_{\left[ 0,s\right) }^{\otimes k}\right) d\tilde{N}_{s}\text{.}
\end{eqnarray*}%
By linearity, for $n\geq 1$, the statement is now completely proved for
every finite linear combination of processes with the form (\ref{prodproc}),
and a standard density argument yields the result for every process with the
form
\begin{equation*}
v_{s}=I_{n}\left( g\left( \cdot ,s\right) \mathbf{1}_{\left[ 0,s\right) \cup
\left( t,1\right] }^{\otimes n}\right) \text{, \ \ }s\in T\text{,}
\end{equation*}%
where the function $g\left( x_{1},...,x_{n},s\right) $ is an element of $%
L^{2}\left( T^{n+1}\right) $ and is symmetric in the variables $\left(
x_{1},...,x_{n}\right) $. To deal with the general case, suppose that $%
u_{s}=\sum_{n\geq 0}I_{n}\left( h_{n}\left( \cdot ,s\right) \right) \in
L^{2}\left( T\times \Omega \right) $ is $\mathbb{F}_{(s,t]^{c}}$ --
predictable on $\left[ 0,t\right] $. This implies, in particular, by setting
$h_{n}^{t}\left( s,t_{1},...,t_{n}\right) $ $=\mathbf{1}_{\left[ 0,t\right]
}\left( s\right) $ $h_{n}\left(t_{1},...,t_{n},s\right) $,%
\begin{eqnarray*}
\sum_{n=0}^{\infty }n!\int_{0}^{t}ds\left\Vert h_{n}\left( \cdot ,s\right)
\right\Vert _{n}^{2} &=&\sum_{n=0}^{\infty }n!\int_{0}^{1}ds\left\Vert
h_{n}^{t}\left( \cdot ,s\right) \right\Vert _{n}^{2}=\sum_{n=0}^{\infty
}n!\left\Vert h_{n}^{t}\right\Vert _{n+1}^{2}<+\infty \text{, \ and} \\
I_{n}\left( h_{n}\left( \cdot ,s\right) \right) &=&I_{n}\left( h_{n}\left(
\cdot ,s \right) \mathbf{1}_{\left[ 0,s\right) \cup \left( t,1\right]
}^{\otimes n}\right) \text{, }n\geq 1\text{,\ }s\in \left[ 0,t\right] \text{.%
}
\end{eqnarray*}%
Now observe that, thanks to the previous discussion,
\begin{equation*}
\left( n+1\right) !\left\Vert \widetilde{h_{n}^{t}}\right\Vert _{n+1}^{2}=%
\mathbf{E}I_{n+1}\left( \widetilde{h_{n}^{t}}\right)
^{2}=n!\int_{0}^{t}\left\Vert h_{n}\left( \cdot ,s\right) \right\Vert
_{n}^{2}ds,
\end{equation*}%
implying
\begin{equation*}
\sum_{n=0}^{\infty }\left( n+1\right) !\left\Vert \widetilde{h_{n}^{t}}%
\right\Vert _{n+1}^{2}=\sum_{n=0}^{\infty }n!\left\Vert h_{n}^{t}\right\Vert
_{n+1}^{2}<+\infty ,
\end{equation*}%
and therefore $u\mathbf{1}_{\left[ 0,t\right] } \in Dom\left( \delta \right)
$. The conclusion is achieved by standard arguments.

(ii) Thanks to \cite[Theorem 4.5]{MPS}, we obtain immediately, for $s\leq t$%
,
\begin{equation*}
G=\mathbf{E}\left( G\mid \mathbb{F}_{s}\right) +\int_{s}^{t}\text{ }^{\left(
p\right) }\left( D_{r}G\right) d\tilde{N}_{r}
\end{equation*}%
where $^{\left( p\right) }\left( \cdot \right) $ indicates the predictable
projection operator with respect to $\mathbb{F}_{s}$, $s\leq t$. To
conclude, it is sufficient to use the independence of the increments of $%
\tilde{N}$ to obtain that, for every $s\leq r\leq t$,%
\begin{eqnarray*}
\mathbf{E}\left( G\mid \mathbb{F}_{s}\right) &=&\mathbf{E}\left( G\mid
\mathbb{F}_{(s,t]^{c}}\right) \text{ \ \ } \\
^{\left( p\right) }\left( D_{r}G\right) &=&\mathbf{E}\left( D_{r}G\mid
\mathbb{F}_{r-}\right) =\mathbf{E}\left( D_{r}G\mid \mathbb{F}%
_{[r,t]^{c}}\right) =\text{ }^{\left( p,t\right) }\left( D_{r}G\right) ,
\end{eqnarray*}%
a.s. -- $\mathbf{P}$. $\ \blacksquare $

\bigskip

\begin{remark}
(i) The arguments used in the proof of Proposition 1-(i) are exclusively
based on the covariance structure of multiple integrals and formula (\ref{r1}%
), and they can be directly applied to the Brownian case. This implies, for
instance, that the Skorohod integral appearing in the statement of
Proposition A.1 in \cite{NP} is also a martingale stochastic integral with
respect to an independently enlarged Brownian filtration.

(ii) We stress that, for the moment, we require the functional $G$, in part
(ii) of Proposition 1, to be $\mathbb{F}_{t}$ - measurable. In the next
section we will show that (\ref{MPS-CO}) holds indeed for every $G\in
\mathbb{D}^{1,2}$. For $t=1$, Proposition 1-(ii) has also been proved in
\cite{Petal}.
\end{remark}

\bigskip

Define $\mathbb{S}^{\ast }$ to be the (dense) subset of $L^{2}\left( T\times
\Omega \right) $ and $\mathbb{L}^{k,2}$, $k\geq 1$, composed of processes of
the type%
\begin{equation*}
v_{s}=\sum_{n=0}^{N}I_{n}\left( f_{n}\left( s,\cdot \right) \right) \text{,
\ \ }s\in T\text{,}
\end{equation*}%
where $N<+\infty $, and, for every $n$, $f_{n}\left(
s,t_{1},...,t_{n}\right) \in S_{n+1}$ and $f_{n}$ is symmetric in the
variables $(t_{1},...,t_{n})$. Then, for every $v\in \mathbb{S}^{\ast }$, a
classic characterization of predictable projections (see \cite[Theorem 6.43]%
{Elliott}) implies immediately, thanks to formulae (\ref{in}) and (\ref{r1}%
), that there exists a jointly measurable application $\phi _{v}$
\begin{equation*}
\Omega \times \Delta ^{2}\mapsto \mathbb{R}:\left( \omega ;s,t\right)
\mapsto \phi _{v}\left( \omega ;s,t\right)
\end{equation*}%
where $\Delta ^{2}=\left\{ \left( s,t\right) \in \left[ 0,1\right]
^{2}:0\leq s\leq t\right\} $, such that, for every $\left( s,t\right) \in
\Delta ^{2}$, $\phi _{v}\left( \cdot ;s,t\right) $ is a version of $\mathbf{E%
}\left[ v_{s}\mid \mathbb{F}_{[s,t]^{c}}\right] $ and, a.s. -- $\mathbf{P}$,
\begin{equation*}
\phi _{v}\left( \cdot ;s,t\right) =\text{ }^{\left( p,t\right) }\left(
v_{s}\right) \text{.}
\end{equation*}

\noindent In general, due again to \cite[Theorem 6.43]{Elliott}, for every
fixed $t\in \left( 0,1\right] $ and every process $u\in L^{2}\left( T\times
\Omega \right) $, the associated predictable projection $^{\left( p,t\right)
}\left( u_{s}\right) $ is such that, for every $s\in \left[ 0,t\right] $, $%
\mathbf{E}\left[ u_{s}\mid \mathbb{F}_{[s,t]^{c}}\right] =$ $^{\left(
p,t\right) }\left( u_{s}\right) $, a.s. -- $\mathbf{P}$. In the future, when
considering the stochastic process $\mathbf{E}\left[ u_{\cdot }\mid \mathbb{F%
}_{[\cdot ,t]^{c}}\right] $, we will implicitly refer to its predictable
modification $^{\left( p,t\right) }\left( u_{\cdot }\right) $. For instance,
with such a convention, formula (\ref{MPS-CO}) can be unambiguously
rewritten as%
\begin{equation*}
G=\mathbf{E}\left( G\mid \mathbb{F}_{(s,t]^{c}}\right) +\int_{s}^{t}\text{ }%
\mathbf{E}\left[ u_{r}\mid \mathbb{F}_{[r,t]^{c}}\right] d\tilde{N}_{r}.
\end{equation*}

\section{Forward-backward martingales and approximation of anticipating
integrals}

In this section, we explore the connection between the anticipating
integrals of the form (\ref{unu}), and a special class of usual It\^{o}
integrals. This relation is applied to prove that, just as on the Wiener
space, anticipating integral processes can be represented as the limit,
under a certain norm, of linear combinations of products of forward and
backward martingales. We start by adapting to the Poisson situation some
known results on the Wiener space. In particular, we will need the following
generalized Clark-Ocone formula, which extends Proposition 1 above and is
the actual equivalent, on the Poisson space, of \cite[Proposition A.1]{NP}.

\begin{prop}
Let the notation of the previous section prevail, and let $G\in \mathbb{D}%
^{1,2}$. Then, for every $0\leq s\leq t\leq 1$, formula (\ref{MPS-CO})
holds. We have also the relation, for $s>0$,%
\begin{equation}
G=\mathbf{E}\left( G\mid \mathbb{F}_{[s,t]^{c}}\right) +\int_{s-}^{t}\text{ }%
^{\left( p,t\right) }\left( D_{r}G\right) d\tilde{N}_{r}=\mathbf{E}\left(
G\mid \mathbb{F}_{[s,t]^{c}}\right) +\delta \left( \text{ }^{\left(
p,t\right) }\left( D_{\cdot }G\right) \mathbf{1}_{\left[ s,t\right] }\left(
\cdot \right) \right)  \label{secondCO}
\end{equation}%
where $\int_{s-}^{t}$ $^{\left( p,t\right) }\left( D_{r}G\right) d\tilde{N}%
_{r}=\lim_{\alpha \uparrow s}\int_{\alpha }^{t}$ $^{\left( p,t\right)
}\left( D_{r}G\right) d\tilde{N}_{r}$.
\end{prop}

\noindent \textbf{Proof. }First observe that the second equality in (\ref%
{secondCO}) follows from the $\mathbb{F}_{(\cdot ,t]^{c}}$ -- predictability
of $^{\left( p,t\right) }\left( D_{\cdot }G\right) $, and an application of
Proposition 1-(i). Moreover, thanks to the martingale property of stochastic
integrals, it is sufficient to prove the statement for $s=0$ and $t\in
\left( 0,1\right] $. We start by considering a random variable $G\in \mathbb{%
D}^{1,2}$ of the form
\begin{equation}
G=I_{m}\left( h\mathbf{1}_{\left[ 0,t\right] }^{\otimes m}\right) \times
I_{n}\left( g\mathbf{1}_{\left( t,1\right] }^{\otimes n}\right)
:=G_{1}\times G_{2}\text{, \ \ }n,m\geq 0\text{,}  \label{simplefunct}
\end{equation}%
where, for $n,m\geq 1$, $h\in L_{s}^{2}\left( \left[ 0,1\right] ^{m}\right) $%
, $g\in L_{s}^{2}\left( \left[ 0,1\right] ^{n}\right) $, and $I_{0}$ stands
for a real constant. Random variables such as (\ref{simplefunct}) are total
in $\mathbb{D}^{1,2}$. Moreover, we can apply Proposition 1 to $G_{1}$ and
obtain, thanks to the stochastic independence between $\mathbb{F}_{t^{c}}$
and $\mathbb{F}_{t}$ and by (\ref{ip})
\begin{eqnarray*}
G &=&G_{2}\times \left[ \mathbf{E}\left( G_{1}\mid \mathbb{F}_{t^{c}}\right)
+\int_{0}^{t}\text{ }^{\left( p,t\right) }\left( D_{r}G_{1}\right) d\tilde{N}%
_{r}\right] \\
&=&\mathbf{E}\left( G\mid \mathbb{F}_{t^{c}}\right) +\int_{0}^{t}\text{ }%
G_{2}\times \text{\ }^{\left( p,t\right) }\left( D_{r}G_{1}\right) d\tilde{N}%
_{r} \\
&=&\mathbf{E}\left( G\mid \mathbb{F}_{t^{c}}\right) +\int_{0}^{t}\text{\ }%
^{\left( p,t\right) }\left( G_{2}\times D_{r}G_{1}\right) d\tilde{N}_{r}%
\text{.}
\end{eqnarray*}%
Note that the last equality comes from \cite[Corollary 6.44]{Elliott} and
from the fact that $G_{2}$ is $\mathbb{F}_{t^{c}}$ measurable, implying that
the (constant) application $r\mapsto G_{2}$, $r\in \left[ 0,t\right] $, can
be interpreted as a $\mathbb{F}_{(r,t]^{c}}$ - predictable process. Finally,
we observe that (\ref{chain rule}) and (\ref{r1}) imply that
\begin{equation*}
G_{2}\times D_{r}G_{1}=D_{r}G\text{, \ \ for every }r\in \left[ 0,t\right] ,
\end{equation*}%
so that, by linearity, (\ref{MPS-CO}) is completely proved for every finite
linear combination of random variables such as (\ref{simplefunct}). Now
suppose that a certain sequence $G^{\left( n\right) }\in \mathbb{D}^{1,2}$
enjoys property (\ref{MPS-CO}) and that $G^{\left( n\right) }$ converges to $%
G$ in $\mathbb{D}^{1,2}$ as $n$ goes to $+\infty $. Then, $\mathbf{E}\left(
G^{\left( n\right) }\mid \mathbb{F}_{t^{c}}\right) \rightarrow \mathbb{E}%
\left( G\mid \mathbb{F}_{t^{c}}\right) $ in $L^{2}\left( \mathbf{P}\right) $%
, and moreover the relation
\begin{equation*}
^{\left( p,t\right) }\left( D_{r}H\right) =\mathbf{E}\left[ D_{r}H\mid
\mathbb{F}_{\left[ r,t\right] ^{c}}\right] \text{,}
\end{equation*}%
true for every fixed $r\in \left[ 0,t\right] $ and every $H\in \mathbb{D}%
^{1,2}$, implies immediately, thanks to Jensen inequality and the isometric
properties of It\^{o} integrals,%
\begin{eqnarray*}
\mathbf{E}\left\{ \left[ \int_{0}^{t}\text{ }\left[ ^{\left( p,t\right)
}\left( D_{r}G^{\left( n\right) }\right) -\text{ }^{\left( p,t\right)
}\left( D_{r}G\right) \right] d\tilde{N}_{r}\right] ^{2}\right\}
&=&\int_{0}^{t}\mathbf{E}\left[ \text{ }\left[ ^{\left( p,t\right) }\left(
D_{r}\left( G^{\left( n\right) }-G\right) \right) \right] ^{2}\right] dr \\
&\leq &\mathbf{E}\int_{0}^{t}\left( D_{r}\left( G^{\left( n\right)
}-G\right) \right) ^{2}dr\rightarrow 0,
\end{eqnarray*}%
and therefore
\begin{equation*}
G=\mathbf{E}\left( G\mid \mathbb{F}_{t^{c}}\right) +\int_{0}^{t}\text{ }%
^{\left( p,t\right) }\left( D_{r}G\right) d\tilde{N}_{r}.
\end{equation*}%
To obtain (\ref{MPS-CO}), use the totality in $\mathbb{D}^{1,2}$ of random
variables such as (\ref{simplefunct}). Eventually, to prove (\ref{secondCO})
just write, for $\varepsilon >0$,%
\begin{equation*}
G=\mathbf{E}\left( G\mid \mathbb{F}_{(s-\varepsilon ,t]^{c}}\right)
+\int_{s-\varepsilon }^{t}\text{ }^{\left( p,t\right) }\left( D_{r}G\right) d%
\tilde{N}_{r}
\end{equation*}%
so that, by letting $\varepsilon \downarrow 0$, the conclusion follows from
the fact that the paths of It\^{o} stochastic integrals (with respect to $%
\tilde{N}$) are c\`{a}dl\`{a}g, as well as relation (\ref{Vfields}) and a
standard martingale argument. $\blacksquare $

\bigskip

\begin{remark}
Proposition 2 can also be proved along the same lines of the proof of \cite[%
Proposition A.1]{NP}. Suppose indeed that $G$ admits the chaotic
decomposition
\begin{equation*}
G=\sum_{m\geq 0}I_{m}(g_{m})\text{, \ \ }g_{m}\in L_{s}^{2}\left(
T^{m}\right) \text{.}
\end{equation*}%
Then, by (\ref{deriv}) and (\ref{r1}), $D_{r}G=\sum_{m\geq
1}mI_{m-1}(g_{m}(\cdot ,r))$ and
\begin{equation*}
\mathbf{E}\left( D_{r}G\mid \mathbb{F}_{[r,t]^{c}}\right) =\sum_{m\geq
1}mI_{m-1}\left( g_{m}(\cdot ,r)\mathbf{1}_{[r,t]^{c}}^{\otimes
(m-1)}\right) .
\end{equation*}%
Since the symmetrization in the $m$ variables $r,t_{1},\ldots ,t_{m-1}$ of
the function $\mathbf{1}_{[s,t]}(r)\mathbf{1}_{[r,t]^{c}}^{\otimes
(m-1)}(t_{1},\ldots ,t_{m-1})$ is given by $\frac{1}{m}\mathbf{1}_{A_{m}}$
where $A_{m}=\bigcup_{i=1}^{m}\{(t_{1},\ldots ,t_{m})\in T^{m},$ $t_{i}\in
[s,t]\}$ we get
\begin{eqnarray*}
\delta \left( \mathbf{1}_{[s,t]}(\cdot )\mathbf{E}\left( D_{\cdot }G\mid
\mathbb{F}_{[\cdot ,t]^{c}}\right) \right) &=&\sum_{m\geq 1}I_{m}(g_{m}%
\mathbf{1}_{A_{m}}) \\
&=&\sum_{m\geq 0}I_{m}(g_{m})-\sum_{m\geq 0}I_{m}(g_{m}\mathbf{1}%
_{A_{m}^{c}}) \\
&=&G-\mathbf{E}\left( G\mid \mathbb{F}_{[s,t]^{c}}\right) .
\end{eqnarray*}
The discussion of Paragraph 2.2 can be used to interpret the Skorohod
integral on the left side as an It\^o integral of a predictable process.
\end{remark}

\vskip0.5cm

The next Proposition shows that every indefinite anticipative integral $%
\int_{0}^{t}u_{s}\delta \tilde{N}_{s}$ can be written, at \textit{fixed} $%
t\in \lbrack 0,1]$, as an \emph{It\^{o}-Skorohod integral} with the form $%
\int_{0}^{t}\mathbf{E}\left( w_{s}\mid \mathbb{F}_{[s,t]^{c}}\right) d\tilde{%
N}_{s}$ where $w$ can be explicitly given in terms of $u$.

\begin{prop}
Let $X$ be a Skorohod integral process $X_{t}=\delta ({u}1_{[0,t]})$, $t\in %
\left[ 0,1\right] $, with $u\in \mathbb{L}^{k,2}$, $k\geq 3$. Then, there
exists a unique process $w\in \mathbb{L}^{k-2,2}$, independent of $t$, such
that, for every fixed $t$,
\begin{equation}
X_{t}=\int_{0}^{t}\mathbf{E}\left( w_{s}\mid \mathbb{F}_{[s,t]^{c}}\right) d%
\tilde{N}_{s}=\delta \left( \mathbf{E}\left( w_{\cdot }\mid \mathbb{F}%
_{[\cdot ,t]^{c}}\right) \mathbf{1}_{\left[ 0,t\right] }\left( \cdot \right)
\right), \mathrm{{a.s. - } \mathbf{P}.}  \label{ItoSkorep}
\end{equation}
\end{prop}

\textbf{Proof. }By applying the Clark-Ocone type formula (\ref{secondCO}) to
the integrand $u$ we can write

\begin{eqnarray*}
X_{t} &=&\int_{0}^{t}u_{s}d\tilde{N}_{s} \\
&=&\int_{0}^{t}\mathbf{E}\left( u_{s}\mid \mathbb{F}_{[s,t]^{c}}\right) d%
\tilde{N}_{s}+\int_{0}^{t}\left( \int_{s-}^{t}\mathbf{E}\left(
D_{r}u_{s}\mid \mathbb{F}_{[r,t]^{c}}\right) \delta \tilde{N}_{r}\right)
\delta \tilde{N}_{s}.
\end{eqnarray*}%
Using a Fubini type theorem (that we can argue exactly as in \cite{NZ}, by
using working on the chaotic expansions) we can interchange the two Skorohod
integrals appearing in the second term to obtain that
\begin{eqnarray*}
X_{t} &=&\int_{0}^{t}\mathbf{E}\left( u_{s}\mid \mathbb{F}%
_{[s,t]^{c}}\right) d\tilde{N}_{s}+\int_{0}^{t}\left( \int_{0}^{r}\mathbf{E}%
\left( D_{r}u_{s}\mid \mathbb{F}_{[r,t]^{c}}\right) \delta \tilde{N}%
_{s}\right) \delta \tilde{N}_{r} \\
&=&\int_{0}^{t}\mathbf{E}\left( u_{s}\mid \mathbb{F}_{[s,t]^{c}}\right) d%
\tilde{N}_{s}+\int_{0}^{t}\mathbf{E}\left( \int_{0}^{r}D_{r}u_{s}\delta
\tilde{N}_{s}\mid \mathbb{F}_{[r,t]^{c}}\right) \delta \tilde{N}_{r} \\
&=&\int_{0}^{t}\mathbf{E}\left( w_{s}\mid \mathbb{F}_{[s,t]^{c}}\right) d%
\tilde{N}_{s}
\end{eqnarray*}%
where we used Proposition 1-(i) as well as the fact that the increments of
the Poisson process are independent and we adopted the notation
\begin{equation}
w_{s}:=u_{s}+\delta (D_{s}u_{\cdot }\mathbf{1}_{[0,s]}(\cdot )):=u_{s}+v_{s}.
\label{w}
\end{equation}%
Let us show that the process $w$ introduced in (\ref{w}) belongs to $\mathbb{%
L}^{k-2,2}$. It suffices indeed to prove that $v\in \mathbb{L}^{k-2,2}$.
Thanks to the inequality (\ref{meyer}) and formula (\ref{comm}), we can
write
\begin{eqnarray*}
\Vert v\Vert _{1,2}^{2} &=&\mathbf{E}\int_{0}^{1}\delta (D_{s}u\mathbf{1}%
_{[0,s]})^{2}ds+\mathbf{E}\int_{0}^{1}\int_{0}^{1}\left( D_{\alpha }\delta
(D_{s}u\mathbf{1}_{[0,s]})\right) ^{2}dsd\alpha \\
&\leq &2\mathbf{E}\int_{0}^{1}\int_{0}^{1}(D_{r}u_{s})^{2}drds+\mathbf{E}%
\int_{0}^{1}\int_{0}^{1}\int_{0}^{1}(D_{\alpha }D_{r}u_{s})^{2}dsdrd\alpha \\
&&+\mathbf{E}\int_{0}^{1}\int_{0}^{1}\delta (D_{\alpha }D_{s}u\mathbf{1}%
_{[0,s]})^{2}dsd\alpha \\
&\leq &2\mathbf{E}\int_{0}^{1}\int_{0}^{1}(D_{r}u_{s})^{2}drds+2\mathbf{E}%
\int_{0}^{1}\int_{0}^{1}\int_{0}^{1}(D_{\alpha }D_{r}u_{s})^{2}dsdrd\alpha \\
&&+\mathbf{E}\int_{0}^{1}\int_{0}^{1}\int_{0}^{1}\int_{0}^{1}(D_{\beta
}D_{\alpha }D_{s}u_{r})^{2}drdsd\alpha d\beta <+\infty ,
\end{eqnarray*}%
since $u\in \mathbb{L}^{3,2}$. In general, it can be similarly proved that
\begin{equation*}
\Vert v\Vert _{k-2,2}^{2}\leq C_{k}\Vert v\Vert _{k,2}^{2},\hskip0,7cmk\geq
3,
\end{equation*}%
where $C_{k}$ is a positive constant depending exclusively on $k$.
Concerning the uniqueness, let us suppose that there exists another process $%
w^{\prime }\in \mathbb{L}^{k-2,2}$ such that $X_{t}=\int_{0}^{t}\mathbf{E}%
\left( w_{s}^{\prime }\mid \mathbb{F}_{[s,t]^{c}}\right) d\tilde{N}_{s}$.
Then, if $z_{s}=w_{s}-w_{s}^{\prime }$, we get
\begin{equation*}
\int_{0}^{t}\mathbf{E}\left( z_{s}\mid \mathbb{F}_{[s,t]^{c}}\right) d\tilde{%
N}_{s}=0
\end{equation*}%
and therefore, for every $t\in (0,1]$,
\begin{equation*}
\int_{0}^{t}\mathbf{E[E}\left( z_{s}\mid \mathbb{F}_{[s,t]^{c}}\right)
]^{2}ds=0.
\end{equation*}%
Let us assume that $z$ has the chaotic expression $z_{s}=\sum_{m\geq
0}I_{m}(g_{m}(\cdot ,s))$ where $g_{m}\in L^{2}(T^{m+1})$; note that, for $%
m\geq 2$,\ the function $g_{m}\left( s_{1},...,s_{m},s\right) $ can be taken
to be symmetric in the variables $\left( s_{1},...,s_{m}\right) $. The above
condition ensures that
\begin{equation*}
\sum_{m\geq 0}\int_{0}^{t}\left( \int_{\left[ s,t\right] ^{c}}\ldots \int_{%
\left[ s,t\right] ^{c}}g_{m}^{2}(s_{1},\ldots ,s_{m},s)ds_{1}\ldots
ds_{m}\right) ds=0
\end{equation*}%
and thus
\begin{equation*}
\int_{\left[ s,t\right] ^{c}}\ldots \int_{\left[ s,t\right]
^{c}}g_{m}^{2}(s_{1},\ldots ,s_{m},s)ds_{1}\ldots ds_{m}=0
\end{equation*}%
for every $m\geq 0$ and for almost every $s\in \lbrack 0,t]$. By letting $%
s\rightarrow t$ we get that $g_{m}\left( \cdot \right) =0$ almost everywhere
on $T^{m+1}$ and the conclusion is obtained. \hfill \vrule width.25cm
height.25cm depth0cm\smallskip

\vskip0.5cm

\begin{remark}
(i) We claim that it is possible to prove a converse to Proposition 3. More
precisely, as in \cite{Tudor}, we can show that if, for a fixed $t$, $%
Y_{t}=\int_{0}^{t}\mathbf{E}\left( w_{s}\mid \mathbb{F}_{[s,t]^{c}}\right) d%
\tilde{N}_{s}$ with $w\in \mathbb{L}^{1,2}$, then there exists $u\in
Dom(\delta )$ such that $Y_{t}=\int_{0}^{t}u_{s}\delta \tilde{N}_{s}$. The
proof would use an analogue of the characterization of Skorohod integrals
stated in Proposition 2.1. of \cite{DN} for the Wiener case; one can prove a
similar characterization in the Poisson context, by following the same line
of reasoning as in \cite[Proposition 2.1]{DN}. This point will be discussed
in a separate paper.

(ii) A question that is not likely to be answered as easily as in the Wiener
case (see \cite{Tudor}), is when a Skorohod integral $X_{t}$ as in (\ref{unu}%
) and an It\^{o}-Skorohod integral such as $Y_{t}=\int_{0}^{t}\mathbf{E}%
\left( w_{s}\mid \mathbb{F}_{[s,t]^{c}}\right) d\tilde{N}_{s}$ are \textit{%
indistinguishable} as stochastic processes, and not only reciprocal
modifications. This is not immediate to answer since, on the Poisson space,
we do not know any sufficiently general criterion, ensuring that an
anticipating integral admits a (right)-continuous version. Nevertheless, as
shown in the following example, one can sometimes apply classical results
from the general theory of stochastic processes.
\end{remark}

\textbf{Example} (Indefinite integrals that are indistinguishable from It%
\^{o}-Skorohod processes). Let the process $X_{t}$, $t\in[0,1]$, be defined
as in (\ref{unu}), and assume that the integrand $u$ and all its Malliavin
derivatives $D^{k}u$ are bounded by a deterministic constant, uniformly on $%
T^{k}\times \Omega$. Then, the assumptions of Proposition 3 are verified,
and we immediately deduce the existence of an It\^{o} -Skorohod integral
with the form $Y_{t}=\int_{0}^{t}\mathbf{E}\left( w_{s}\mid \mathbb{F}%
_{[s,t]^{c}}\right) d\tilde{N}_{s}$, $t\in[0,1]$, such that $X$ and $Y$ are
modifications. We claim that, in this setting, $X$ and $Y$ also admit two
indistinguishable c\`adl\`ag modifications. Recall indeed the following
classic criterion (see \cite[Ch. III, Section 4]{GiSk}): if $Z$ is a
stochastic process such that for some $p\geq 1$ and $\beta >0$,
\begin{equation}
\mathbf{E}\left| (Z_{t+h}-Z_{t})(Z_{t}-Z_{t-h})\right| ^{p} \leq C(p)\times
h^{1+\beta } ,  \label{cadlagC}
\end{equation}
where $C(p)$ is a positive constant, then $Z$ admits a c\`adl\`ag
modification. Note that if $Z^{\prime}$ is a modification of $Z$, and $Z$
satisfies (\ref{cadlagC}), then (\ref{cadlagC}) must also hold for $%
Z^{\prime}$. We shall prove (\ref{cadlagC}) for $Z=X$. To this end, note $%
I_{h+} = [t,t+h]$ and $I_{h-} = [t-h,t]$, so that, under the above
assumptions and by (\ref{dua}),
\begin{eqnarray*}
&&\mathbf{E}\left| (X_{t+h}-X_{t})(X_{t}-X_{t-h})\right| ^{2} = \mathbf{E}%
\left( \delta (u\mathbf{1}_{I_{h+}}) \delta (u\mathbf{1}_{I_{h+}})\delta (u%
\mathbf{1}_{I_{h-}})^{2}\right) \\
&&= \mathbf{E}\int _{I_{h+}}u_{s} D_{s} \left[ \delta (u\mathbf{1}%
_{I_{h+}})\delta (u\mathbf{1}_{I_{h-}})^{2}\right] ds = \mathbf{E}\left(
\int _{I_{h+}}u_{s} \left[ D_{s} \delta (u\mathbf{1}_{I_{h+}})\right] ds
\right)\delta (u\mathbf{1}_{I_{h-}})^{2} \\
&& +\mathbf{E}\left( \int _{I_{h+}}u_{s} \left[ D_{s}\delta (u\mathbf{1}%
_{I_{h-}})^{2}\right] ds\right) \delta (u\mathbf{1}_{I_{h+}}) +\mathbf{E}%
\left( \int _{I_{h+}}u_{s} \left[D_{s} \delta (u\mathbf{1}_{I_{h+}})\right] %
\left[ D_{s}\delta (u\mathbf{1}_{I_{h-}})^{2}\right] ds \right) \\
&&:= A+B+C.
\end{eqnarray*}
By formula (\ref{comm}), the first summand above can be decomposed as
follows
\begin{eqnarray*}
A&=&\mathbf{E}\int _{I_{h+}}u_{s}\left[ u_{s}\mathbf{1}_{h+}(s) + \delta (
D_{s}u \mathbf{1}_{h+}) \right] ds := A_{1}+ A_{2}.
\end{eqnarray*}
We need only show how to handle $A_{1}$, and similar techniques can be used
to deal with the remaining terms $A_{2},B$ and $C$. Indeed, we can write
\begin{equation*}
A_{1} =\mathbf{E}\delta (u\mathbf{1}_{I_{h-}})^{2}\int _{I_{h+}}u^{2}_{s}ds
\leq cst\times h \mathbf{E}\delta (u\mathbf{1}_{I_{h-}})^{2}\leq cst\times
h^{2}
\end{equation*}
since, by (\ref{meyer}),
\begin{equation*}
\mathbf{E}\delta (u\mathbf{1}_{I_{h-}})^{2}\leq \mathbf{E}\int
_{I_{h-}}u^{2}_{s}ds+ \int _{I_{h-}}\int _{I_{h-}}(D_{s}u_{\alpha } )
^{2}dsd\alpha \leq cst\times h
\end{equation*}
because $u$ and its derivative are assumed to be uniformly bounded. Since (%
\ref{cadlagC}) is also true for $Z=Y$, we deduce that there exist two
processes $X^{\prime}$ and $Y^{\prime}$ such that $X^{\prime}$ is a
c\`adl\`ag modification of $X$ and $Y^{\prime}$ is a c\`adl\`ag modification
of $Y$. We can now apply a classic argument (see for instance \cite[Theorem
2, p. 4]{Protter}), to deduce that $X^{\prime}$ and $Y^{\prime}$ are also
indistinguishable. $\blacksquare$

\vskip0.5cm

Let us recall some notation taken from \cite{PTT}. By $L_{0}^{2}(\mathbf{P})$
we denote the set of zero mean square integrable random variables. We will
write $\mathbf{BF}$ for the class of stochastic processes that can be
expressed as finite linear combinations of processes of the type
\begin{equation*}
Z_{t}=\mathbf{E}\left( H_{1}\mid \mathbb{F}_{t}\right) \times \mathbf{E}%
\left( H_{2}\mid \mathbb{F}_{t^{c}}\right) =M_{t}\times M_{t}^{\prime }
\end{equation*}%
where $H_{1}\in L_{0}^{2}(\mathbf{P})$ and $H_{2}\in L^{2}(\mathbf{P})$.
Plainly, $M$ is a martingale with respect to $\mathbb{F}_{t}$ and $M^{\prime
}$ is a \textit{backward martingale}. By backward martingale we mean that $%
M_{t}^{\prime }$ is in $L^{1}\left( \mathbf{P}\right) $ and $\mathbb{F}%
_{t^{c}}$ -- measurable for every $t$, and $\mathbf{E}\left( M_{s}^{\prime
}\mid \mathbb{F}_{t^{c}}\right) =M_{t}^{\prime }$ for any $s\leq t$; see
e.g. \cite{RY}.

\smallskip

We give a counterpart of Lemma 2 in \cite{PTT}. The proof needs a slightly
different argument.

\begin{lemma}
Let $A_{1},A_{2}$ be two disjoint Borel subsets of $[0,1]$ and assume that $%
F $ is a random variable in $\mathbb{D}^{k,2}$, $k\geq 1$, such that $F$ is
measurable with respect to the $\sigma $-algebra $\mathbb{F}_{A_{1}}\vee
\mathbb{F}_{A_{2}}$. Then, $F$ is the limit in $\mathbb{D}^{k,2}$ of linear
combinations of smooth random variables of the type%
\begin{equation}
G=G_{1}\times G_{2}\text{,}  \label{prod}
\end{equation}%
where, for $i=1,2$, $G_{i}$ is a polynomial, $\mathbb{F}_{A_{i}}$ -
measurable functional.
\end{lemma}

\textbf{Proof.} Suppose first that $F$ is a simple functional of the form
\begin{equation}
F=\tilde{N}(h_{1})\ldots \tilde{N}(h_{n})  \label{s1}
\end{equation}%
where $n\geq 1$ and $h_{i}\in L^{2}(T)$, $i=1,...,n$. Then, the conclusion
can be obtained exactly as in Lemma 2 of \cite{PTT}, by using twice formula (%
\ref{r1}). The next step is to consider $F=I_{n}(f)$ with $f\in L^{2}(T^{n})$%
. In this case, by the above discussed definition of the multiple integral $%
I_{n}$, $F$ is the limit in $\mathbb{D}^{k,2}$ of random variables $F_{k}$
as in (\ref{s1}). Let us denote by $p_{m,k}$ a sequence of linear
combinations of product (\ref{prod}) such that $p_{m,k}\rightarrow F_{k}$ in
$\mathbb{D}^{k,2}$ as $m\rightarrow \infty $. Clearly,
\begin{equation*}
\Vert F-p_{m,k}\Vert _{k,2}\leq \Vert F-F_{k}\Vert _{k,2}+\Vert
F_{k}-p_{m,k}\Vert _{k,2}
\end{equation*}%
and this goes to zero when $m,k\rightarrow \infty $. Eventually, take the
general case $F=\sum_{m\geq 0}I_{m}(f_{m})$ where $f_{m}\in L^{2}(T^{n})$
are symmetric functions. The conclusion will follow if we prove that $F$ can
be approximated in $\mathbb{D}^{k,2}$ by a sequence $F^{N}$ of random
variables with finite chaotic expansion and this is trivial if we put $%
F^{N}=\sum_{n=0,...,N}I_{n}(f_{n}).$ See also Proposition 1.2.1 in \cite{N}
for further details. \hfill \vrule width.25cm height.25cm depth0cm\smallskip

\smallskip

\begin{remark}
Recall the relation between multiple integrals and Charlier polynomials
stated in formula (\ref{CharlierP}). Then, by inspection of the proof of
Lemma 1 (which partially follows that of Lemma 2 in \cite{PTT}), and thanks
to the multiplication formula (\ref{product}), it is clear that if $%
F=I_{n}(f)$ ($\in \mathbb{D}^{k,2}$), then $F$ can be approximated in $%
\mathbb{D}^{k,2}$ by linear combinations of random variables with the form $%
C_{n}\left( t_{k+1}-t_{k},\tilde{N}_{t_{k+1}}-\tilde{N}_{t_{k}}\right) $
where $0\leq \ t_{k}<t_{k+1}\leq 1$, and $C_{n}$ is the $n$th Charlier
polynomial. Therefore, the random variables $G_{i}$, $i=1,2$, appearing in
the proof of Lemma 1 can be chosen as polynomial functionals of degrees $%
d_{i}$, $i=1,2$, such that $d_{1}+d_{2}\leq n$.
\end{remark}

\vskip0.5cm

We shall also introduce the following quadratic variation (in mean) of a
given measurable process $\left\{ X_{t}:t\in T\right\} $ such that $\mathbf{E%
}X_{t}^{2}<+\infty $ for every $t$:%
\begin{equation*}
V(X)=\sup_{\pi }\mathbf{E}\sum_{i=0}^{n-1}\left(
X_{t_{i+1}}-X_{t_{i}}\right) ^{2},
\end{equation*}%
where $\pi $ runs over all partitions of $T=\left[ 0,1\right] $, with the
form $0=t_{0}<t_{1}<\ldots <t_{n}=1$.

\vskip0.5cm

We state the main result of this section.

\begin{theorem}
Let $X$ be a Skorohod integral process $X_{t}=\delta ({u}\mathbf{1}_{[0,t]})$
with $u\in \mathbb{L}^{k,2}$, $k\geq 3$. Then there exists a sequence of
processes $(Z_{t}^{(r)})_{t\in \lbrack 0,1]}$, $r\geq 1$ such that
\begin{equation}
Z^{(r)}\in \mathbf{BF}\mbox{ for every }r\geq 1  \label{t1}
\end{equation}%
and
\begin{equation}
\lim_{r\rightarrow \infty }V\left( X-Z^{(r)}\right) =0.  \label{t2}
\end{equation}
\end{theorem}

\textbf{Proof.} Remark first that $V(X)<\infty $ by Proposition 1 in \cite%
{TV}. We will use the It\^{o}-Skorohod representation of $X_{t}=\int_{0}^{t}%
\mathbf{E}\left( w_{s}\mid \mathbb{F}_{[s,t]^{c}}\right) d\tilde{N}_{s}$
with $w$ given by (\ref{w}). For $n\geq 1$ and a partition $\pi =\{0=t_{0}$ $%
<...<t_{n}$ $=1\}$, we introduce the approximation of $w$%
\begin{equation}
w_{t}^{\pi }=\sum_{i=0}^{n-1}\frac{1}{t_{i+1}-t_{i}}\left(
\int_{t_{i}}^{t_{i+1}}\mathbf{E}\left( w_{s}\mid \mathbb{F}%
_{[t_{i},t_{i+1}]^{c}}\right) ds\right)
1_{(t_{i},t_{i+1}]}(t):=\sum_{i=0}^{n-1}F_{i}\mathbf{1}_{(t_{i},t_{i+1}]}(t),%
\hskip0.5cmt\in \lbrack 0,1].  \label{wpi}
\end{equation}%
Since $w\in \mathbb{L}^{1,2}$ then $w^{\pi }\in \mathbb{L}^{1,2}$ and $%
w^{\pi }$ converges to $w$ in $\mathbb{L}^{1,2}$ as $|\pi |\rightarrow 0$
(see \cite{TV}, and also \cite{N} for the Gaussian case). Note that the
random variables $F_{i}$, $i=0,...,n$, appearing in (\ref{wpi}) are
measurable with respect to $\mathbb{F}_{[t_{i},t_{i+1}]^{c}}$. We also set,
for $\pi $ as above,%
\begin{equation*}
Y_{t}^{\pi }=\delta \left( \mathbf{E}\left( w_{\cdot }^{\pi }\mid \mathbb{F}%
_{[\cdot ,t]^{c}}\right) \mathbf{1}_{\left[ 0,t\right] }\right) =\int_{0}^{t}%
\mathbf{E}\left( w_{s}^{\pi }\mid \mathbb{F}_{[s,t]^{c}}\right) d\tilde{N}%
_{s}\text{, \ \ }t\in T.
\end{equation*}%
Using properties (\ref{r1}), (\ref{r2}) and (\ref{ip}) of Poisson Skorohod
integrals we therefore deduce%
\begin{eqnarray}
Y_{t}^{\pi } &=&\sum_{i=0}^{n-1}\int_{0}^{t}1_{(t_{i},t_{i+1}]}(s)\mathbf{E}%
\left( F_{i}\mid \mathbb{F}_{[s,t]^{c}}\right) \delta \tilde{N}_{s}
\label{relazione} \\
&=&\sum_{i=0}^{n-1}\int_{0}^{t}1_{(t_{i},t_{i+1}]}(s)\mathbf{E}\left(
F_{i}\mid \mathbb{F}_{[t_{i},t_{i+1}\vee t]^{c}}\right) \delta \tilde{N}_{s}
\notag \\
&=&\sum_{i=0}^{n-1}\mathbf{E}\left( F_{i}\mid \mathbb{F}_{(t_{i},t_{i+1}\vee
t]^{c}}\right) \left( \tilde{N}_{t\wedge t_{i+1}}-\tilde{N}_{t_{i}}\right)
\mathbf{1}_{(t\geq t_{i})};  \notag
\end{eqnarray}%
note, in particular, that the last equality in (\ref{relazione}) derives
from an application of formula (\ref{ip}), where the last two terms vanish
thanks to (\ref{r2}) (alternatively, one can also use Proposition 1-(i)).
This is all we need to conclude the proof of Theorem 1. As a matter of fact,
as in \cite[proof of Theorem 1]{PTT}, we can now use Lemma 1 to prove that $%
Y^{\pi }$ (and hence $X$) can be approximated, in the sense of formula (\ref%
{t2}), by a sequence of processes $Z^{(r)}$ satisfying (\ref{t1}). \hfill %
\vrule width.25cm height.25cm depth0cm\smallskip

\vskip0.5cm

We state a converse result to Theorem 1; the arguments of \cite[Theorem 4]%
{PTT} apply, and the proof is therefore omitted. It shows that the "$V$%
-norm" is somewhat complete.

\begin{theorem}
Let $Z^{(n)}\in \mathbf{BF}$, $n\geq 1$ be such that $V(Z^{(n)})<\infty $
and $Z^{(n)}$ is a Cauchy sequence in the $V$-norm, in the sense that
\begin{equation*}
\lim_{n,m\rightarrow \infty }V\left( Z^{(n)}-Z^{(m)}\right) =0.
\end{equation*}%
Then, there exists a Skorohod integral process $X$ with $V(X)<\infty $ such
that
\begin{equation*}
\lim_{n\rightarrow \infty }V\left( Z^{(n)}-X\right) =0.
\end{equation*}
\end{theorem}

\section{On the stochastic calculus for anticipating integrals on the
Poisson space}

In the previous Section we have seen that, for every fixed $t\in \lbrack 0,1]
$, the Skorohod integral $X_{t}=\int_{0}^{t}u_{s}\delta \tilde{N}_{s}$ is
equal to the It\^{o}-Skorohod integral $Y_{t}=\int_{0}^{t}\mathbf{E}\left(
w_{s}\mid \mathbb{F}_{[s,t]^{c}}\right) d\tilde{N}_{s}$. As seen in Section
2, the last random variable is the It\^{o} integral of the predictable
process $\mathbf{E}\left( w_{\cdot }\mid \mathbb{F}_{[\cdot ,t]^{c}}\right) $%
, with respect to the $\mathbb{F}_{(\cdot ,t]^{c}}$ -- martingale $\tilde{N}%
_{\cdot }$. As such, it is an isometry, in the sense that
\begin{equation*}
\mathbf{E}\left( \int_{0}^{t}\mathbf{E}\left( w_{s}\mid \mathbb{F}%
_{[s,t]^{c}}\right) d\tilde{N}_{s}\right) ^{2}=\mathbf{E}\int_{0}^{t}\mathbf{%
E}\left( w_{s}\mid \mathbb{F}_{[s,t]^{c}}\right) ^{2}ds
\end{equation*}%
(this can be also derived from formula (\ref{cov})), and can moreover be
approximated by a sequence of martingales. Define indeed
\begin{equation}
Y_{t}^{\lambda }:=\int_{0}^{\lambda }\mathbf{E}\left( w_{s}\mid \mathbb{F}%
_{[s,t]^{c}}\right) d\tilde{N}_{s}=\delta \left( \mathbf{1}_{[0,\lambda
]}\left( \cdot \right) \mathbf{E}\left( w.\mid \mathbb{F}_{[\cdot
,t]^{c}}\right) \right) .  \label{ylambdat}
\end{equation}%
Then, for every fixed $t\in \lbrack 0,1]$, the process $\left\{
Y_{t}^{\lambda }:0\leq \lambda \leq t\right\} $ is a martingale
with respect to the filtration $\mathbb{F}_{(\lambda ,t]^{c}}$,
$\lambda \leq t$, and it holds that $Y_{t}^{t}=Y_{t}$, and for any
$\lambda <t$,
\begin{eqnarray*}
\mathbf{E}\left\vert Y_{t}^{\lambda }-Y_{t}\right\vert ^{2} &=&\mathbf{E}%
\left\vert \delta \left( \mathbf{1}_{(\lambda ,t]}\left( \cdot \right)
\mathbf{E}\left( w.\mid \mathbb{F}_{[\cdot ,t]^{c}}\right) \right)
\right\vert ^{2} \\
&=&\mathbf{E}\int_{\lambda }^{t}\mathbf{E}\left( w_{s}\mid \mathbb{F}%
_{[s,t]^{c}}\right) ^{2}ds
\end{eqnarray*}%
which goes to $0$ as $\lambda \rightarrow t$ by the dominated convergence
theorem. As a consequence, the $Y_{t}^{\lambda }$ converges in $L^{2}(%
\mathbf{P})$ to $Y_{t}$ as $\lambda \rightarrow t$ and, by a standard
martingale convergence theorem (see e.g. Problem 3.20 in \cite{KS}), the
convergence holds a.s. -- $\mathbf{P}$. This fact allows us to introduce a
stochastic calculus of It\^{o} type for the It\^{o}-Skorohod integral $Y_{t}$
(and hence for indefinite Skorohod integrals $X_{t}$). The main idea is to
use the tools of the stochastic calculus for the martingale $Y_{t}^{\lambda }
$ and to let $\lambda \rightarrow t$. We obtain in this way a change of
variable formula for the indefinite integral processes; this seems quite
interesting since, as far as we know, there is no It\^{o} formula \emph{\`{a}
la Nualart-Pardoux} \cite{NP} for anticipating integrals in the Poisson
case. We also derive a Burkholder-type bound for the $L^{p}$-norm of a
Skorohod integral.

\vskip0.5cm

\begin{prop}
(It\^{o}'s formula ) Let $f\in C^{2}(\mathbb{R})$ , fix $t\in T$, and define
$Y_{t}=\delta \left( \mathbf{E}\left( w_{\cdot }\mid \mathbb{F}_{[\cdot
,t]^{c}}\right) \mathbf{1}_{\left[ 0,t\right] }\left( \cdot \right) \right) $%
, where $w\in L^{2}(\Omega \times T)$. Then it holds that
\begin{eqnarray}
f(Y_{t}) &=&f(0)+\int_{0}^{t}f^{\prime }(Y_{t}^{s-})\mathbf{E}\left(
w_{s}\mid \mathbb{F}_{[s,t]^{c}}\right) d\tilde{N}_{s}  \label{ItoFormule} \\
&&+\frac{1}{2}\int_{0}^{t}f^{\prime \prime }(Y_{t}^{s-})\mathbf{E}\left(
w_{s}\mid \mathbb{F}_{[s,t]^{c}}\right) ^{2}ds  \notag \\
&&+\sum_{0\leq s\leq t}\left( f(Y_{t}^{s})-f(Y_{t}^{s-})-f^{\prime
}(Y_{t}^{s-})\left( Y_{t}^{s}-Y_{t}^{s-}\right) \right)  \notag
\end{eqnarray}%
where $Y_{t}^{s-}=\lim_{\alpha \rightarrow s,\alpha <s}Y_{t}^{\alpha }.$ In
particular, let $X_{t}=\int_{0}^{t}u_{s}\delta \tilde{N}_{s}$, $t\in T$,
where $u\in \mathbb{L}^{k,2}$, $k\geq 3$, and let $w$ be the process
appearing in Proposition 3, formula (\ref{ItoSkorep}); then, for every $t\in
T$, $f\left( X_{t}\right) $ equals the right-hand side of (\ref{ItoFormule}).
\end{prop}

\textbf{Proof. }Fix $t\in \left( 0,1\right] $, and consider the
process
\begin{equation*}
\widetilde{Y}_{t}^{\lambda }=\left\{
\begin{array}{ll}
Y_{t}^{\lambda } & \text{if }\lambda \leq t \\
Y_{t} & \text{if }t<\lambda <+\infty%
\end{array}%
\right.
\end{equation*}%
as well as the family of $\sigma $-fields%
\begin{equation*}
\widetilde{\mathbb{F}}_{\lambda }^{t}=\left\{
\begin{array}{ll}
\mathbb{F}_{(\lambda ,t]^{c}} & \text{if }\lambda \leq t \\
\mathbb{F}_{1} & \text{if }t<\lambda <+\infty%
\end{array}%
\right.
\end{equation*}%
then, the application $\lambda \rightarrow \widetilde{Y}_{t}^{\lambda }$
defines a square integrable c\`{a}dl\`{a}g martingale $\widetilde{\mathbb{F}}%
_{\lambda }^{t}$ with respect to the filtration $\widetilde{\mathbb{F}}%
_{\lambda }^{t}$, $\lambda \geq 0$. Moreover,
\begin{eqnarray*}
\widetilde{Y}_{t}^{\lambda } &=&\int_{0}^{\lambda }\mathbf{E}\left(
w_{s}\mid \mathbb{F}_{[s,t]^{c}}\right) \mathbf{1}_{\left( s\leq t\right) }d%
\tilde{N}_{s}\text{ and} \\
\text{ }\left\langle \widetilde{Y}_{t}^{\cdot },\widetilde{Y}_{t}^{\cdot
}\right\rangle _{\lambda } &=&\int_{0}^{t\wedge \lambda }\mathbf{E}\left(
w_{s}\mid \mathbb{F}_{[s,t]^{c}}\right) ^{2}ds.
\end{eqnarray*}%
We can therefore apply It\^{o}'s formula (see e.g. \cite{Protter}, Theorem
32, p. 71) at $\lambda =t$ to obtain%
\begin{eqnarray*}
f(\widetilde{Y}_{t}^{t}) &=&f\left( Y_{t}\right) =f(0)+\int_{0}^{t}f^{\prime
}(\widetilde{Y}_{t}^{s-})\mathbf{E}\left( w_{s}\mid \mathbb{F}%
_{[s,t]^{c}}\right) d\tilde{N}_{s} \\
&&+\int_{0}^{t}f^{\prime \prime }(\widetilde{Y}_{t}^{s-})\mathbf{E}\left(
w_{s}\mid \mathbb{F}_{[s,t]^{c}}\right) ^{2}ds \\
&&+\sum_{0<s\leq t}\left( f(\widetilde{Y}_{t}^{s})-f(\widetilde{Y}%
_{t}^{s-})-f^{\prime }(\widetilde{Y}_{t}^{s-})\left( \widetilde{Y}_{t}^{s}-%
\widetilde{Y}_{t}^{s-}\right) \right) . \\
&=&f(0)+\int_{0}^{\lambda }f^{\prime }(Y_{t}^{s-})\mathbf{E}\left( w_{s}\mid
\mathbb{F}_{[s,t]^{c}}\right) d\tilde{N}_{s} \\
&&+\frac{1}{2}\int_{0}^{t}f^{\prime \prime }(Y_{t}^{s-})\mathbf{E}\left(
w_{s}\mid \mathbb{F}_{[s,t]^{c}}\right) ^{2}ds \\
&&+\sum_{0<s\leq t}\left( f(Y_{t}^{s})-f(Y_{t}^{s-})-f^{\prime
}(Y_{t}^{s-})\left( Y_{t}^{s}-Y_{t}^{s-}\right) \right) ,
\end{eqnarray*}%
by the definition of $\widetilde{Y}$.

\hfill \vrule width.25cm height.25cm depth0cm\smallskip

\bigskip

\vskip0.5cm

We now show that Proposition 4 can be applied to write a change of variables
formula for stochastic processes that are representable as the product of a
martingale and a backward martingale.

\begin{prop}
Let $M$(resp. $M^{\prime }$) be a martingale (resp. a backward martingale)
with respect to the filtration $\left\{ \mathbb{F}_{t}:t\in T\right\} $, and
suppose moreover that $M_{1}\in \mathbb{D}^{1,2}$, $\mathbf{E}\left(
M_{0}\right) =0$ and%
\begin{equation*}
\mathbf{E}\left( \int \left( a_{s}M_{s}^{\prime }\right) ^{2}ds\right)
<+\infty ,
\end{equation*}
where $a_{r}=\mathbf{E}\left( D_{r}M_{1}\mid \mathbb{F}_{r-}\right) .$ Then,
for every $f\in C^{2}(\mathbb{R})$ we have
\begin{eqnarray*}
f(M_{t}M_{t}^{\prime }) &=&f(0)+\int_{0}^{t}f^{\prime }(M_{s-}M_{t}^{\prime
})M_{t}^{\prime }a_{s}d\tilde{N}_{s} \\
&&+\frac{1}{2}\int_{0}^{t}f^{\prime \prime }(M_{s-}M_{t}^{\prime
})(M_{t}^{\prime })^{2}a_{s}^{2}ds \\
&&+\sum_{0\leq s\leq t}\left[ f(M_{s}M_{t}^{\prime })-f(M_{s-}M_{t}^{\prime
})-f^{\prime }(M_{s-}M_{t}^{\prime })M_{t}^{\prime }(M_{s}-M_{s-})\right] .
\end{eqnarray*}
\end{prop}

\textbf{Proof.} Remark that a product $M_{t}M_{t}^{\prime }$ is a It\^{o}%
-Skorohod integral. Indeed, by a standard Clark-Ocone formula (see \cite[%
Theorem 4.5]{MPS})
\begin{equation*}
M_{t}=\mathbf{E}\left( M_{1}\mid \mathbb{F}_{t}\right) =\int_{0}^{t}\mathbf{E%
}\left( D_{r}M_{1}\mid \mathbb{F}_{r-}\right) d\tilde{N}_{r}
\end{equation*}%
and using (\ref{ip}) and (\ref{r2}), it holds that
\begin{eqnarray*}
M_{t}M_{t}^{\prime } &=&M_{t}^{\prime }\int_{0}^{t}a_{s}d\tilde{N}%
_{s}=\int_{0}^{t}a_{s}M_{t}^{\prime }d\tilde{N}_{s} \\
&=&\int_{0}^{t}a_{s}\mathbf{E}\left( M_{s}^{\prime }\mid \mathbb{F}%
_{[s,t]^{c}}\right) d\tilde{N}_{s}=\int_{0}^{t}\mathbf{E}\left(
a_{s}M_{s}^{\prime }\mid \mathbb{F}_{[s,t]^{c}}\right) d\tilde{N}_{s}
\end{eqnarray*}%
By applying Proposition 4, we obtain
\begin{eqnarray*}
M_{t}M_{t}^{\prime } &=&f(0)+\int_{0}^{t}f^{\prime }(Z_{t}^{s-})\mathbf{E}%
\left( a_{s}M_{s}^{\prime }\mid \mathbb{F}_{[s,t]^{c}}\right) d\tilde{N}_{s}
\\
&&+\frac{1}{2}\int_{0}^{t}f^{\prime \prime }(Z_{t}^{s-})\mathbf{E}\left(
a_{s}M_{s}^{\prime }\mid \mathbb{F}_{[s,t]^{c}}\right) ^{2}ds \\
&&+\sum_{0\leq s\leq t}\left( f(Z_{t}^{s})-f(Z_{t}^{s-})-f^{\prime
}(Y_{t}^{s-})\left( Z_{t}^{s}-Z_{t}^{s-}\right) \right)
\end{eqnarray*}%
where $Z_{t}^{\lambda }=\int_{0}^{\lambda }\mathbf{E}\left(
a_{s}M_{s}^{\prime }\mid \mathbb{F}_{[s,t]^{c}}\right) d\tilde{N}_{s}$. The
conclusion follows, since for every $\lambda \leq t$, we can write

\begin{eqnarray*}
Z_{t}^{\lambda } &=&\mathbb{E}\left( \int_{0}^{\lambda }\mathbf{E}\left(
a_{s}M_{s}^{\prime }\mid \mathbb{F}_{[s,\lambda ]^{c}}\right) d\tilde{N}%
_{s}\mid \mathbb{F}_{[\lambda ,t]^{c}}\right) \\
&=&\mathbb{E}\left( M_{\lambda }M_{\lambda }^{\prime }\mid \mathbb{F}%
_{[\lambda ,t]^{c}}\right) =M_{\lambda -}M_{t}^{\prime }.
\end{eqnarray*}%
and
\begin{equation*}
Z_{t}^{s-}=\lim_{\alpha \uparrow s}Z_{t}^{\alpha -}=M_{t}^{\prime }M_{s-}.
\end{equation*}%
\hfill \vrule width.25cm height.25cm depth0cm\smallskip

\vskip0.5cm

Here is a more particular situation.

\begin{corollary}
If $f\in C^{2}(\mathbb{R})$, we have for every $t\in T$
\begin{eqnarray*}
f(\tilde{N}_{t}(\tilde{N}_{1}-\tilde{N}_{t})) &=&f(0)+(\tilde{N}_{1}-\tilde{N%
}_{t})\int_{0+}^{t}f^{\prime }\left( \tilde{N}_{s-}(\tilde{N}_{1}-\tilde{N}%
_{t})\right) d\tilde{N}_{s} \\
&&+\frac{1}{2}(\tilde{N}_{1}-\tilde{N}_{t})^{2}\int_{0+}^{t}f^{\prime
}\left( \tilde{N}_{s-}(\tilde{N}_{1}-\tilde{N}_{t})\right) ds \\
&&+\sum_{0\leq s\leq t\leq 1}\left[ f\left( \tilde{N}_{s}(\tilde{N}_{1}-%
\tilde{N}_{t})\right) -f\left( \tilde{N}_{s-}(\tilde{N}_{1}-\tilde{N}%
_{t})\right) \right.  \\
&&\left. -f^{\prime }\left( \tilde{N}_{s-}(\tilde{N}_{1}-\tilde{N}%
_{t})\right) (\tilde{N}_{1}-\tilde{N}_{t})(\tilde{N}_{s}-\tilde{N}_{s-})%
\right] .
\end{eqnarray*}
\end{corollary}

\textbf{Proof.} Apply Proposition 5 with $M_{t}=\tilde{N}_{t}$, $%
M_{t}^{\prime }=\tilde{N}_{1}-\tilde{N}_{t}$ and $a\equiv 1$. \hfill \vrule %
width.25cm height.25cm depth0cm\smallskip

\vskip0.5cm

We conclude this section by proving a class of Burkholder type inequalities.
These could an useful tool to bound the $L^{p}$-norms of anticipating
integrals, since on the Poisson space there are no analogous of Meyer's
inequalities (see \cite{N}) proved for the operators $D$ and $\delta $ as
defined in Section 2.

\begin{prop}[Burkholder inequalities]
If $Y_{t}=\int_{0}^{t}\mathbf{E}\left( w_{s}\mid \mathbb{F}%
_{[s,t]^{c}}\right) d\tilde{N}_{s}$ with $w\in L^{2}(T\times \Omega
)$, then, for every $p\geq 1$ there exist two universal constants
$K_{1}(p)>0$ and $K_{2}(p)>0$ such that
\begin{equation*}
K_{1}(p)\mathbf{E}\left( \int_{0}^{t}\left( \mathbf{E}\left( w_{s}\mid
\mathbb{F}_{[s,t]^{c}}\right) \right) ^{2}d[\tilde{N}]_{s}\right) ^{\frac{p}{%
2}}\leq \mathbb{E}\left\vert Y_{t}\right\vert ^{p}\leq K_{2}(p)\mathbf{E}%
\left( \int_{0}^{t}\left( \mathbf{E}\left( w_{s}\mid \mathbb{F}%
_{[s,t]^{c}}\right) \right) ^{2}d[\tilde{N}]_{s}\right) ^{\frac{p}{2}}
\end{equation*}%
where $[\tilde{N}]_{t}=N_{t}$. In particular, let $X_{t}=\int_{0}^{t}u_{s}%
\delta \tilde{N}_{s}$, $t\in T$, where $u\in \mathbb{L}^{k,2}$, $k\geq 3$,
and let $w$ be the process defined in formula (\ref{ItoSkorep}); then, for
every $t\in T$,
\begin{equation*}
K_{1}(p)\mathbf{E}\left( \int_{0}^{t}\left( \mathbf{E}\left( w_{s}\mid
\mathbb{F}_{[s,t]^{c}}\right) \right) ^{2}d[\tilde{N}]_{s}\right) ^{\frac{p}{%
2}}\leq \mathbf{E}\left\vert X_{t}\right\vert ^{p}\leq K_{2}(p)\mathbf{E}%
\left( \int_{0}^{t}\left( \mathbf{E}\left( w_{s}\mid \mathbb{F}%
_{[s,t]^{c}}\right) \right) ^{2}d[\tilde{N}]_{s}\right) ^{\frac{p}{2}}.
\end{equation*}
\end{prop}

\textbf{Proof. } We have, by classical Burkholder inequalities for jump
processes (see e.g. \cite[Theorem 54]{Protter})
\begin{eqnarray*}
\mathbf{E}\left\vert Y_{t}\right\vert ^{p} &\leq &\mathbf{E}\sup_{\lambda
\leq t}\left\vert Y_{t}^{\lambda }\right\vert ^{p} \\
&\leq &K_{2}(p)\mathbf{E}\left( \int_{0}^{\lambda }\left( \mathbf{E}\left(
w_{s}\mid \mathbb{F}_{[s,t]^{c}}\right) \right) ^{2}d[\tilde{N}]_{s}\right)
^{\frac{p}{2}} \\
&\leq &K_{2}(p)\mathbf{E}\left( \int_{0}^{t}\left( \mathbf{E}\left(
u_{s}\mid \mathbb{F}_{[s,t]^{c}}\right) \right) ^{2}d[\tilde{N}]_{s}\right)
^{\frac{p}{2}}.
\end{eqnarray*}%
For the lower bound, we write
\begin{eqnarray*}
&&\mathbf{E}\left\vert Y_{t}\right\vert ^{p}=\mathbf{E}\lim_{\lambda
\rightarrow t}\left\vert Y_{t}^{\lambda }\right\vert ^{p}=\lim_{\lambda
\rightarrow t}\mathbf{E}\left\vert Y_{t}^{\lambda }\right\vert ^{p} \\
&\geq &\lim_{\lambda \rightarrow t}K_{1}(p)\mathbf{E}\left(
\int_{0}^{\lambda }\left( \mathbf{E}\left( w_{s}\mid \mathbb{F}%
_{[s,t]^{c}}\right) \right) ^{2}d[\tilde{N}]_{s}\right) ^{\frac{p}{2}}=%
\mathbf{E}\left( \int_{0}^{t}\left( \mathbf{E}\left( w_{s}\mid \mathbb{F}%
_{[s,t]^{c}}\right) \right) ^{2}d[\tilde{N}]_{s}\right) ^{\frac{p}{2}}.
\end{eqnarray*}%
\hfill \vrule width.25cm height.25cm depth0cm

\end{document}